\documentclass[psamsfonts]{amsart}
\usepackage{amsfonts}
\usepackage{amsmath,amsxtra,amssymb,latexsym,amscd ,amsthm}
\usepackage{ucs}
\usepackage[mathscr]{eucal}
\usepackage{graphicx}
\usepackage{color}
\usepackage{graphics}
\usepackage{xypic,tikz}
\usepackage[all]{xy}
\UseComputerModernTips
\usepackage{hyperref}
\hypersetup{
    colorlinks= false,
    linkcolor=blue,
    filecolor=magenta,      
    urlcolor=cyan,
}
 
\usepackage{enumitem}
\usepackage{tikz-cd}

\usepackage{fancyhdr}

\usepackage{setspace}
\numberwithin{equation}{section}

\newtheorem{thm}{Theorem}[section]
\newtheorem{lem}[thm]{Lemma}

\newtheorem{dfn}[thm]{Definition}
\newtheorem{cor}[thm]{Corollary}
\newtheorem{exa}[thm]{Example}

\theoremstyle{remark}
\newtheorem{rem}[thm]{Remark}

\title{Tangent cones and $C^1$ regularity of definable sets}

\author{Krzysztof Kurdyka, Olivier Le Gal and Nhan Nguyen}

\address{LAMA, Universit\'e Savoie Mont Blanc, 73376 Le Bourget-du-Lac Cedex, France}
\email{Krzysztof.Kurdyka@univ-savoie.fr}
\email{Olivier.Le Gal@univ-savoie.fr}
\email{nguyenxuanvietnhan@gmail.com }

\newcommand{\bb}{\mathbb}
\newcommand{\al}{\mathcal}
\newcommand{\bff}{\mathbf}

\newcommand{\sspan}{{\rm span}}

\newcommand{\ptg}{{\rm ptg}}
\newcommand{\tg}{{\rm tg}}
\newcommand{\gr}{{\rm \Gamma}}
\newcommand{\card}{{\rm card}}

\begin{document}
\maketitle

\parskip .12cm
\begin{abstract}
Let $X\subset \bb R^n$ be a connected locally closed definable set in an o-minimal structure. We prove that the following three statements are equivalent: (i) $X$ is a $C^1$ manifold, (ii) the tangent cone and the paratangent cone of $X$ coincide at every point in $X$, (iii) for every $x \in X$, the tangent cone of $X$ at the point $x$ is a $k$-dimensional linear subspace of $\bb R^n$ ($k$ does not depend on $x$) varies continuously in $x$, and the density $\theta(X, x) < 3/2$.
\end{abstract}

\section{Introduction}
Let $X$ be a subset of $\bb R^n$ and let $x \in \bb R^n$. The \textit{tangent cone $\tg_xX$  and paratangent cone $\ptg_x X$ of $X$ at the point $x$} are defined as follows: if $x \not \in \overline{X}$, $\tg_xX = \ptg_x X = \{0\}$, and otherwise, 
$$\tg_x X := \{ a u \mid a \in \bb R, a \geq 0,  u = \lim_{i\to \infty} \frac{x_i - x}{\|x_i - x\|}, \{x_i\}\subset X, \{x_i\} \to x\},$$
$${\rm ptg}_x X : = \{a u\mid a \geq 0,   u = \lim_{i \to \infty} \frac{x_i - y_i}{\|x_i - y_i\|}, X \supset \{x_i\} \to x, X\supset \{y_i\} \to x\}.$$
Note that $\tg_x X$ and $\ptg_x X$ are closed sets in $\bb R^n$.
We denote by $\tg X:= \{(x, v), x \in X, v \in \tg_x X\}$ and  $\ptg X := \{(x,  v), x \in X, v\in \ptg_x X\}$. 

Characterizing $C^1$ submanifolds of $\bb R^n$ in terms of their tangent cones has been studied by many authors, see for example \cite{gluck}, \cite{tierno}, \cite{bi-gr}, \cite{gh-ho}, or a survey of  Bigolin and Golo \cite{bi-go}. In this paper we restrict ourselves to this problem in the context of o-minimal structures. We first prove that a connected locally closed definable subset of $\bb R^n$ is a $C^1$ manifold if and only if its tangent cone and paratangent cone coincide at every point (Theorem \ref{thm_two_cones_definable}). This result is a strong version of the two-cones coincidence theorem (Theorem \ref{thm_two_cones}) which was initially proved by Tierno \cite{tierno}. The result is no longer true if definability is omitted (Remark \ref{rem_definable_two_cones}). 

Next, we discuss a result recently established by Ghomi and Howard \cite{gh-ho} that if $X\subset \bb R^n$ is a locally closed set such that for each $x \in X$, $\tg_xX$ is a hyperplane (i.e., a $(n-1)$ linear subspace of $\bb R^n$), and varies continuously in $x$, then $X$ is a union of $C^1$ hypersurfaces. Moreover, if the lower density $\Theta(X, x)$ is at most $m < \frac{3}{2}$ for every $x \in X$ then $X$ is a $C^1$ hypersurface. A natural question thus arises here is whether the result  remains true if $\tg_x X$ are $k$-planes with $k < n-1$. 

In section 4, we show in Example \ref{exa_main} that the first statement in the result of Ghomi-Howard is not always true if $k < (n-1)$. We also prove that the second statement is still valid, more precisely that if $X$ is a locally closed definable set such that for every $x\in X$, $\tg_x X$ is a $k$-dimensional linear subspace ($k$ is independent of $x$) varying continuously in $x$ and the density $\theta(X, x) < 3/2$ (need not be upper bounded by an $m < \frac{3}{2})$, then $X$ is a $C^1$ manifold (Theorem \ref{thm_gh_ho_general}).  Notice that, in general, notions of lower density $\Theta(X, x)$  and density $\theta(X,x )$ are different. Nevertheless, with the conditions on tangent cones as above, they coincide. Therefore, our result can be considered as a generalization of the result of Ghomi-Howard.

Throughout the paper, $\bb R^n$ denotes the $n$-dimensional Euclidean space equipped with the standard norm $\|x\| = \sqrt{x_1^2 + \ldots + x_n^2}$ where $x = (x_1, \ldots, x_n)$; $\overline{\bff B^n}(x, r)$, $\bff B^n(x, r)$ and  $\bff S^{n-1}(x, r)$) denote respectively the closed ball, the open ball and the sphere in $\bb R^n$ of radius $r$ centered at $x$.  Let $X$ be a subset in $\bb R^n$. Denote by $\overline{X}$ the closure of $X$ in $\bb R^n$ and by $\partial X : = \overline{X}\setminus X$ the boundary of $X$. Let $f$ be a map. We denote by $\gr_f$ the graph of $f$.

The Grassmanian $\bb G_n^k$ of all $k$-dimensional linear subspaces of $\bb R^n$ is endowed with the metric $\delta$ defined as follows: for $P, Q$ in $\bb G_n^k$,
$$\delta(P, Q) = \sup_{v \in P, \|v\| = 1}\{ \|v - \pi_Q (v)\|\},$$
where $\pi_Q: \bb R^n \to Q$ is the orthogonal projection from $\bb R^n$ to $Q$.
Following \cite{gh-ho}, we say that $P\in \bb G_n^k$ and $Q\in \bb G_n^k$ are \emph{orthogonal} when $\delta(P,Q)=1$. Remark that this terminology does not coincide with the usual notion of orthogonality for general subspaces in Euclidean geometry :  $P =\{(x,y,z);\; x=0\}$ and $Q=\{(x,y,z);\; y = 0\}$ are orthogonal according to our definition but not all vectors in $P$ are orthogonal to any vector in $Q$.  

By a \textit{$k$-dimensional $C^1$ manifold  in $\bb R^n$} (or $C^1$ manifold for simplicity) we mean a subset of $\bb R^n$, locally $C^1$ diffeomorphic to $\bb R^k$; \textit{a hypersurface} in $\bb R^n$  is a $C^1$ manifold in $\bb R^n$ of dimension $n - 1$. 

Let $X \subset \bb R^n$. In the paper, we often denote by $\pi_x : \bb R^n \to \tg_x X$ the orthogonal projection. By abuse of notation, here  we identify $\tg_x X$ with its translation $\{x + \tg_x X\}$.

By a \textit{definable set} we mean a set which is definable (with parameters) in an o-minimal expansion $(\bb R , <, +, \cdot, \dots )$ of the ordered field of real numbers. Definable sets form a large class of subsets of $\mathbb{R}^n$: for instance, any semi-algebraic set, any sub-analytic set is definable. We refer the reader to \cite{vdd}, \cite{coste} for the basic properties of o-minimal structures. In the paper, we will use Curve selection Lemma (\cite{coste}, Theorem 3.2), Uniform finiteness on fibers (\cite{coste} Theorem 2.9) and Hardt's definable triviality Theorem (\cite{coste}, Theorem 5.22) without repeating the references.

\section{Bundle of vector spaces}

Let $X$ be a subset of $\bb R^n$. Let $E \subset X \times \bb R^n$. For $x \in X$ we denote by 
$$E_x: =\{v\in \bb R^n: (x, v) \in E\}$$ 
the \textit{fiber} of $E$ at the point $x$. For $U\subset X$, we set
$$E|_U: = \{(x, v) \in E: x \in U\}$$ and call it \textit{the restriction of $E$ to $U$}.
If every fiber of $E$ is a linear subspace of $\bb R^n$ we call $E$ a \textit{bundle of vector spaces over $X$}, or a \textit{bunlde over $X$}, or just a \textit{bundle} if the base $X$ is clear from the context. We call $E$ a\textit{ trivial bundle} if all its fibers have the same dimension, and a \textit{closed bundle} if it is a closed set in $\bb R^{2n}$. 

Suppose $E$ is a trivial bundle over $X$. If the map $X \to \bb G_n^k$ defined by $x \mapsto E_x$ is continuous, we say that $E$ is a \textit{continuous trivial bundle}.


\begin{lem}[\cite{glaeser}, Propositions I, II, pages 39, 40]\label{lem_glaeser_0}

Let $E\subset X \times \bb R^n$.
\begin{enumerate}

\item If $E$ is a closed trivial bundle then $E$ is continuous.

\item If $E$ is a closed bundle then the function $x \mapsto \dim E_x$ is upper-semicontinuous, i.e., for $x \in X$ there is an open neighborhood $U_x$ of $x$ in $X$ such that $\dim E_x \geq \dim E_y$ for every $y \in U_x$.
\end{enumerate}

\end{lem}

\begin{proof}
$(i)$ - Suppose that $k$ is the dimension of fibers of $E$. If $E$ is not continuous, there exists a sequence $\{x_i\} \subset X$ tending to $x$, $\lim_{i\to \infty}  E_{x_i} \to \tau \in \bb G_n^k$ and $\tau \not\subset E_x$. By the closedness of $E$, $\tau \subset E_x$, which is a contradiction.

$(ii)$ - Suppose the assertion is not true, i.e., there exist a point $x \in X$ and a sequence $\{x_m\}$ in $X$ tending to $x$ such that $\dim E_{x_m} > \dim E_x$. We may assume that $\lim_{m \to \infty} E_{x_m}  = P \in \bb G_n^k$, since $\bb G_n^k$ is compact. Note that $k > \dim E_x$. Since $E$ is closed,  $P \subset E_x$. This implies $k \leq \dim E_x$, which is a contradiction.   
\end{proof}

\begin{lem}\label{lem_paratangent} Let $X$ be a locally closed subset of $\bb R^n$. If $\ptg X$ is a trivial bundle then $\ptg X$ is continuous. 
\end{lem}

\begin{proof} It follows directly from the definition of the paratangent cone that  $\ptg(\overline{X})$ is a closed set in $\bb R^{2n}$ and  $\ptg_x X = \ptg_x \overline{X}$ for every $x \in X$. If $V \subset X$ is a closed set in $\bb R^n$ then $\ptg X|_{V} = \ptg \overline{X}|_{V}$ is a closed set in $\bb R^{2n}$.

Let $x \in X$. Since $X$ is locally closed, there is $W_x$, a neighborhood of $x$ in $X$, which is closed in $\bb R^n$. The restriction $\ptg X|_{W_x}$ is then a closed set in $\bb R^{2n}$. Since $\ptg X$ is a trivial bundle, so is $\ptg X|_{W_x}$. By  $(i)$ in Lemma \ref{lem_glaeser_0} $\ptg X|_{W_x} $ is continuous,  $\ptg X$ is, therefore, a continuous trivial bundle. 
\end{proof}
 
\section{Two-cones coincidence Theorem}

The aim of this section is to prove Theorem \ref{thm_two_cones_definable}, a strong version of two-cones coincidence theorem of Tierno for definable sets. 

We need the following two lemmas which generalize Lemma 3.3 and Lemma 3.1 in \cite{gh-ho}.

\begin{lem}\label{lem_open_map}
Let $X$ be a locally closed subset of $\bb R^n$ such that $\tg X$ is a continuous trivial bundle of $k$-dimensional vector spaces. Let $x \in X$ and $H$ be a  $k$-plane in $\bb R^n$ which is not orthogonal to $\tg_x X$. Let $\pi: \bb R^n \to H$ the orthogonal projection. Then, there exists an open set $U$ of $x$ in $X$ such that $\pi|_ {U}: U \to H $ is an open map.   
\end{lem}

\begin{proof} The proof follows closely the proof of Lemma 3.3, \cite{gh-ho}.

By the continuity of $\tg X$, we can choose an open neighborhood $U$ of $x$ in $X$ such that for all $q\in U$, $\tg_q X$ is not orthogonal to $H$, or equivalently that $\tg_q X$ is transverse to $H^\perp$, the orthogonal complement of $H$ in $\bb R^n$. We will prove that $\pi|_U$ is an open map. Fix $q \in U$. By the local closedness of $X$, there is an $r > 0$ small enough such that $W: = X \cap \overline{\bff B}^n(q, r)\subset U$ is a compact set. Moreover, the boundary $\partial W$ is in $\partial \overline{\bff B}^n(q, r)$, meaning $q \not\in \partial W$.  With $r$ sufficiently small, we may assume that 
$$\pi(q) \not\in \pi(\partial W)$$
because, otherwise, there exists a sequence of positive numbers $\{r_i\}$ tending to $0$ such that for each $i$, there is a point $q_i \in X \cap \partial\overline{\bff B}^n(q, r_i)$ with $\pi(q_i) = \pi(q)$. This implies that $\sspan (q - q_i) \subset H^\perp$. As $i\to \infty$ we have $q_i \to q$ and the sequence $\sspan (q - q_i)$ (extracting a subsequence if necessary) tends to a line $l \in \tg_x X$. Thus, $l \in H^\perp \cap \tg_q X$; but $\dim (H^\perp \cap \tg_q X) = \dim H^\perp + \dim \tg_q X - n = 0$ since  $ H^\perp$ is transverse to $\tg_q X$, which is a contradiction. 

Since $\pi (\partial W)$ is a compact set, there is $s > 0$ such that 
\begin{equation}\label{fm_boundary}
\overline{\bff B}^k(\pi(q), s)\cap  \pi(\partial W) = \emptyset.
\end{equation} 
It suffices  to show that $\pi(W)$ contains an open neighborhood of $\pi(q)$ in $H$.

Suppose on the contrary that $\overline{\bff B}^k(\pi(q), \varepsilon) \not\subset \pi(W)$, $\forall \varepsilon > 0$. Choose a point $q' \in \overline{\bff B}^k(\pi(q), s/2)\setminus \pi(W)$ and let $s'$ be the distance from $q'$ to $\pi(W)$. Note that $s' \leq s/2$. Since $\pi(W)$ is compact, $\overline{\bff B}^k (q', s') \cap \pi(W) \neq \emptyset$ and $\bff B ^k(q', s') \cap \pi(W)  = \emptyset$. For every $p \in\overline{\bff B}^k(q', s')$,  
$$\|p - \pi(q)\| \leq \|p - q'\| + \|q' - \pi(q)\| \leq  s' + s/2 \leq  s.$$
This means that $\overline{\bff B}^k(q', s') \subset \overline{\bff B}^k (\pi(q), s)$. By (\ref{fm_boundary}),  $\overline{\bff B}^k(q', s') \cap \pi(\partial W) = \emptyset$.  Take $y' \in \overline{\bff B}^k(q', s') \cap \pi(W)$ and $y \in \pi^{-1}(y') \cap W$. Note that $y \not \in \partial W$, so $y$ is an interior point of $W$, and hence $\tg_y W = \tg_y X$ which is a linear subspace. 

Since $\bff B^k(q', s') \cap \pi(W)  = \emptyset$, no point of $W$ is contained in the cylinder $C: = \pi^{-1} (\bff B ^k(q', s'))$. This implies that $\tg_y X \subset \tg_y \partial C$. Both $\tg_y X$ and $H^\perp$ are included in $  \tg_y \partial C$, so $$\dim (\tg_y X \cap  H^\perp) \geq \dim \tg_y X + \dim H^\perp - \dim \tg_y \partial C =   k + (n-k) - (n-1) = 1.$$ This shows that $\tg_y X$ is orthogonal to $H$, which is a contradiction. \end{proof}

\begin{lem}\label{lem_C^1_map} Let $U \subset \bb R^k$ be an open set and $f: U \to \bb R^{n-k}$ be a map. Suppose that $ \gr_f$ is locally closed. If $\tg \gr_f$ is a continuous trivial bundle of $k$-dimensional vector spaces and $\tg_{(x, f(x))}\gr_f$ is not orthogonal to $\bb R^k$ for every $x \in U$, then $f$ is $C^1$. 
\end{lem}
\begin{proof}
We first prove that $f$ is continuous.
Suppose on the contrary  that $f$ is not continuous, meaning that there are $x \in U$ and a sequence $\{x_i\}$ in $U$ tending to $x$ such that $\lim_{i\to \infty} f(x_i) = y \neq f(x)$. 

Since $\gr_f$ and the orthogonal projection $\pi: \bb R^n \to \bb R^k$ satisfy the hypothesis of Lemma \ref{lem_open_map}, there is an open set $V$ of $(x, f(x))$ in $\gr_f$ such that $\pi|_V$ is an open map. Set $W : = \pi (V)$, which is an open  neighborhood of $x$ in $\bb R^k$. Since $\{x_i\}$ tends to $x$, there is $N \in \bb N$ such that $x_i \in  W$ for all $i > N$. This implies that $(x_i, f(x_i)) \in V$ for all $i > N$.

If  $y \neq \infty$, shrinking $V$ so that $(x,y) \not \in \overline{V}$, there  is a neighborhood $K$ of $(x,y)$ in $\bb R^n$ such that $K \cap V = \emptyset$. Since $f(x_i)$ tends to $y$,  we have $(x_i, f(x_i)) \in K$ for all $i > N$ when $N$ is large enough. This shows that $K \cap V \neq \emptyset$, which is a contradiction. 

If $y = \infty$, $(x_i, f(x_i)) \not \in V$ for all $i > N$ when $N$ is large enough. This again gives a contradiction.  

We now show that $f$ is a $C^1$ map. 

Let $\{a_1, \ldots, a_k, b_1, \ldots, b_{n-k}\}$ be the canonical basis of $\bb R^n$. For $x \in U$, consider the function $f^i_x (t): = f(x + ta_i)$. The graph $\gr_{f^i_x}$ of $f^i_x$ is the intersection 
$\gr_{f^i_x} = \gr_f \cap (x,0)+\sspan(a_i,b_1,\dots,b_k)$.
This implies that $$\tg_{(x, f(x))} \gr_{f^i_x} \subset \tg_{(x, f(x))}\gr_f \cap \sspan(a_i,b_1,\dots,b_{n-k})=: l_x.$$ But $l_x$ is a line, because $\tg_{(x, f(x))}\gr_f$ is not orthogonal to $\bb R^k\times\{0\}^{n-k}$. On the other hand, since $f$ is continuous, $\gr_{f^i_x}$ is a continuous curve, so $\tg_{(x, f(x))} \gr_{f^i_x}$ has dimension at least $1$. Then $\tg_{(x, f(x))} \gr_{f^i_x} = l_x$, so $f^i_x$ is differentiable at $t=0$. Thus, $f$ has partial derivatives at any point. 

The bundle $\tg_{(x, f(x))}\gr_f$ is continuous, hence $l_x$, its transverse intersection with $\sspan(a_i,b_1,\dots, b_{n-k})$ is continuous. Therefore  $f$ has continuous partial derivatives on $U$, so $f$ is $C^1$.

\end{proof}
\begin{rem}The statement of Lemma 3.1 \cite{gh-ho} is similar to the statement of Lemma \ref{lem_C^1_map} except the local closedness of $ \gr_f$ is missing. This is a gap because $f$ might not be continuous if $ \gr_f$ is not locally closed. For example, consider the function $f(x) = 0$ if $x$ is a rational number, and $f(x) = 1$ otherwise. The tangent cone to $\gr_f$ is the $x$-axis at any point, hence $\gr_f$ is a continuous trivial bundle, but $f$ is not continuous. 
\end{rem}

\begin{thm}\label{thm_gluck} \label{thm_gluck} A locally closed set  $X \subset \bb R^n$  is a $C^1$ manifold if and only if $\tg X$ is a continuous trivial bundle and the restriction of the map $\pi_x : X \to \tg_x X$ to some neighborhood of $x$ in $X$ is injective.
\end{thm}

\begin{proof}
The necessity is a trivial fact. We now prove the sufficiency. For $x \in X$, by the hypothesis, there exists an open neighborhood $U$ of $x$  such that $\varphi: = \pi_x|_{U}: U \to \tg_x X$ is injective. Moreover, $\pi_x$ is open by Lemma \ref{lem_open_map}. This implies $\varphi: U \to \varphi(U)$ is a homeomorphism. Consider the map $\psi:=\varphi^{-1}: \varphi(U) \to U \subset \bb R^n$. We have $\tg  \gr_\psi = \tg X|_{\gr_\psi}$, which is a continuous trivial bundle. Shrinking $U$ if necessary we may assume that $\tg_y X$ is not orthogonal to $\tg_x X$ for every $y \in U$.  The map $\psi$ then satisfies the conditions of Lemma \ref{lem_C^1_map}, so it is of class $C^1$, meaning that $U$ is a $C^1$ manifold. 
\end{proof}
\begin{rem} Theorem \ref{thm_gluck} is slightly stronger than a similar result proved by Gluck (Theorem 10.1, \cite{gluck}). In the result of Gluck, $X$ is assumed to be a topological manifold instead of a locally closed set as in our statement.  
\end{rem}

\begin{thm}[Two-cones coincidence, Tierno \cite{tierno}]\label{thm_two_cones}
A locally closed subset $X$ of $\bb R^n$ is a $C^1$ manifold if and only if $TX$ and $\ptg X$ coincide, and both are trivial bundles of vector spaces over $X$.
\end{thm}
\begin{proof}
Since $\ptg X$ is a trivial bundle, it is continuous by Lemma \ref{lem_paratangent}. On the other hand, $\tg X = \ptg X$, hence $\tg X$ is a continuous trivial bundle. 

Let $x \in X$. We denote by $\pi_x: \bb R^n \to \tg_x X$ the orthogonal projection. By Theorem \ref{thm_gluck}, it suffices to prove that the map $\pi_x$ is injective on some neighborhood of $x$ in $X$. 

Suppose on contrary that there are sequences of points $\{z_i\}_i$ and $\{z'_i\}_i$ in $X$ converging to $x$ such that $ \pi_x(z_i) = \pi_x(z'_i)$ for all $i$. This implies that $\sspan (z_i - z'_i)$ accumulates to a  line $l \subset \tg_x X ^ \perp$. By the definition, $l \subset \ptg_x X $. Since $\ptg_x X = \tg_x X$, $l \subset \tg_x X \cap \tg_x X^\perp = 0$, a contradiction. 
\end{proof}

\begin{thm}[Definable two-cones coincidence]\label{thm_two_cones_definable} A connected, locally closed definable subset of $\bb R^n$ is a $C^1$ manifold if and only if $\tg X$ and $\ptg X$ coincide.
\end{thm}
\begin{proof}

We just need to show the sufficiency. First we prove that for every $x \in X$,  $\tg_xX$ is a linear subspace of $\bb R^n$, or equivalently that $\tg X$ is a bundle. Fix $x \in X$, we may identify $x$ with the origin $0$. By the hypothesis,  $\tg_0X = \ptg_0 X $ which is symmetric, i.e., if $v \in \tg_0 X$ so is $-v$. It is enough to verify that if $v, w \in \tg_0 X$ then $ v + w \in  \tg_0 X$. Since  $v, -w \in \tg_0X$ and $X$ is a definable set there exist two curves $\gamma, \beta$ in $X$ starting at $0$ such that $v \in \tg_0\gamma$ and $-w \in \tg_0 \beta$ (see Curve selection Lemma). Choose sequences of points $\{x_i\}_i\subset \gamma$ and $\{y_i\}_i\subset \beta$ converging to $0$ such that  $\|x_i\| = \frac{a}{b}\|y_i\|$, where $a: = \|v\|$ and $b: = \|w\|$.  Thus, 
$$v + w = \lim_{i \to \infty} (  a \frac{x_i}{\|x_i\|} - b \frac{y_i}{\|y_i\|} )
 = \lim_{i\to \infty} \frac{a}{\|x_i\|} (x_i - y_i) \in \ptg_0X = \tg_0X.$$

From now on, we set
 $$O_k:=\{z \in X: \dim \tg_x X \leq k \}.$$ 

Fix $k$ and let $z \in O_k$. Take $V$ to be a closed neighborhood of $z$ in $X$. Since $X$ is locally closed, we may assume $V$ to be a closed set in $\bb R^n$ so $\ptg X|_V$ is a closed  bundle. By Lemma \ref{lem_glaeser_0} (ii),  for $y \in V$, the map $y \mapsto \dim \ptg_y X$ is upper-semicontinuous, so is the map $ y \mapsto \dim \tg_y X$, meaning that there exists  $U\subset V$, an open neighborhood of $z$, such that $k = \dim \ptg_z X \geq \dim \ptg_y X$ for all $y$ in $U$. This implies that $U \subset O_k$, hence $O_k$ is an open set in $X$.

Since $X$ is definable,  $\dim \tg_x X \leq \dim X = d$ for every $x \in X$ (see \cite{ku-ra}, Lemma 1.2), hence $O_{d} = X$. Set $\al O := X\setminus O_{d-1}$, which is a closed set of $X$.

We denote by $X_{\rm Sing}$ the set of singular points of $X$, i.e., points at which $X$ fails to be a $C^1$ manifold of dimension $d$. Remark that $X_{\rm Sing}$ is a definable set of dimension less than $d$ (see for instance \cite{vdd}, \cite{coste}). 

Since $O_{d-1} \subset X_{\rm Sing}$, $\dim O_{d-1} <d $. For $x \in \al O$, 
$$ \tg_x X \supset \tg_x \al O \supset \tg_x X \setminus \tg_x  O_{d-1}.$$

Taking the closures of all sets above,
$$ \tg_x X \supset \tg_x \al O \supset \overline{\tg_x X \setminus \tg_x  O_{d-1}}.$$

Because $\tg_x X$ is a linear space of dimension $d$ and $\tg_x  O_{d-1}$ is a linear subspace of $\tg_x X$ of dimension less than $d$, $\overline{\tg_x X \setminus \tg_x  O_{d-1}} = \tg_x X$. So, 
$$ \tg_x X = \tg_x \al O \subset \ptg_x \al O \subset\ptg_x X = \tg_x X.$$

Since $\al O$ is closed in the locally closed set $X$, it is also locally closed. As been shown above, $\tg \al O = \ptg \al O = \ptg X|_{\al O}$ which is a trivial bundle. By Theorem \ref{thm_two_cones}, $\al O$ is a $C^1$ manifold of dimension $d$. 
Next we will prove that $X = \al O$, therefore, it is a $C^1$ manifold. 

Let $x \in \al O$ and $\pi_x: \bb R^n \to \tg_x \al O$ be the orthogonal projection.  It follows from Theorem \ref{thm_gluck} and Lemma \ref{lem_open_map} that there is an open neighborhood $U$ of $x$ in $\bb R^n$ such that the restriction of $\pi_x$ to $U \cap \al O$ is injective and $V: = \pi_x (U\cap \al O)$ is an open set in $\tg_x \al O$. Set $ W: = \pi_x^{-1}(V) \cap U$, so that $W$ is an open neighborhood of $x$ in $\bb R^n$ with $\pi_x (W) = V$. This implies that 

$$ \pi_x (W \cap X) = \pi_x (W \cap \al O) = V.$$
Since $\tg_x \al O = \ptg_x X$, shrinking $U$ if necessary, the restriction of $\pi_x $ to $W \cap X$ is injective. The sets $W \cap X$ and $W \cap \al O$ then are graphs of mappings over the same domain $V$. On the other hand, $ W \cap \al O \subset W \cap X$, then $ W \cap \al O = W \cap X$. This means that $W \cap \al O$, an open neighborhood of $x$ in $\al O$, is an open neighborhood of $x$ in $X$. Thus, $\al O$ is an open set in $X$. Since $\al O$ is both closed and open in $X$ and $X$ is connected, $\al O$ is equal to $X$.

\end{proof}

A direct consequence of Theorem \ref{thm_two_cones} is:
\begin{cor} Let $X \subset \bb R^n$ be a locally closed definable set. Suppose that $ \tg_x X = \ptg_x X$ for every $x \in X$. Then, each connected component of $X$ is a $C^1$ manifold.
\end{cor}

\begin{rem}\label{rem_definable_two_cones}
The definability in the hypothesis of Theorem \ref{thm_two_cones_definable} is necessary. Consider the following locally closed sets, 
$$X:= \overline{\{ (x, \sin\frac{1}{x}), x \neq 0\}} \setminus \{(0,1), (0, -1)\} \subset \bb R^2,$$ and
$$Y := \{ -1 -\frac{1}{n}, n \in \bb N \} \cup [-1,1] \cup \{ 1 + \frac{1}{n}, n \in \bb N\}.$$
The sets $X$ and $Y$ are not definable in any o-minimal structure : $X\cap\mathbb R \times \{0\}$ and $Y$ have infinitely many componant. The set $X$ is connected, $\tg X = \ptg X$, but $X$ is not a $C^1$ manifold. The set $Y$ has $\tg Y= \ptg Y$, but $[-1,1]$, a connected component of  $Y$, is not a $C^1$ manifold.
\end{rem}

\section{Definable sets with  continuous trivial tangent cones}
Let us recall the result proved by Ghomi and Howard \cite{gh-ho}.

\begin{dfn}\label{def_lower_density}\rm Let $X\subset \bb R^n$, $x \in \bb R^n$ . Suppose that \textit{ the Hausdorff dimension} of $\tg_x X$, denoted by $\dim_{\al H} \tg_x X$, is an integer $k$. The \textit{lower density} $\Theta(X, x)$ of $X$ at the point $x$ is defined as follows: if $x \not \in \overline{X}$ then $\Theta(X, x) = 0$, and otherwise, 
$$\Theta(X, x) = \lim_{r\to 0} \inf \frac{\al H^k(X \cap \overline{\bff B}^n(x, r))}{r^k \mu_k}$$
where $\al H^k$ is the $k$-dimensional Hausdorff measure, $\mu_k$ is the volume of the unit ball of dimension $k$.
\end{dfn}

\begin{thm}[Theorem 1.1, \cite{gh-ho}]\label{thm_gh_ho} Let $X$ be a locally closed subset of $\bb R^n$. Suppose that $\tg X$ is a $(n-1)$-dimensional continuous trivial bundle. Then,
\begin{enumerate}
\item $X$ is a union of $C^1$ hypersurfaces;
\item if $\Theta(X, x)$ is at most $m < 3/2$ then  $X$ is a $C^1$ hypersurface. 
\end{enumerate}
\end{thm}

The following example shows that in general the statement $(i)$ of Theorem \ref{thm_gh_ho} is no longer true when the hyperplanes are replaced by $k$-planes with $k < n-1$, meaning that a locally closed subset in $\bb R^n$ with continuous trivial tangent cone might not be a union of $C^1$ manifolds. 

\begin{exa}\label{exa_main} \rm
We identify $\bb C$ with $\bb R^2$. Consider the map $h: \bb C \to \bb C$ defined as follows: 
$$h(z) : = \frac{z^2}{|z|^\frac{3}{2}} \text{ if } z \neq 0, \text{ and } h(0) = 0.
$$
Denote by $X\subset \mathbb R^4$ the graph of $h$. 
We have

(1) $X$ is locally closed;

(2) $\tg X$ is a $2$-dimensional  continuous trivial bundle;

(3) $X$ is not a union of $C^1$ submanifolds of dimension $2$ of $\bb R^4$. 

\begin{proof}
Remark that $h$ is continuous, hence $X$ is a topological manifold and the condition (1) is automatically satisfied. Moreover, $h$ is smooth except at the origin, where statement (2) is obvious. We now calculate $\tg_0X$. Let $x \in X \setminus \{0\}$, we may write $x = (z, h(z))$ for some $z \in \bb C^*$. We have
$$
\frac{x}{\|x\|} = \frac{(z, h(z))}{\|(z, h(z))\|}= \left(\frac{z}{\sqrt{|z|}} \frac{1}{\sqrt{1 + |z|}}, \frac{z^2}{|z^2|} \frac{1}{\sqrt{1 + |z|}} \right).
$$
Notice that $\frac{z}{\sqrt{|z|}} \frac{1}{\sqrt{1 + |z|}} \to 0$ when $z \to 0$. Hence $\tg_0X \subset \{z = 0\}$. On the other hand, if $z = re^{i\theta }$ and $r \to 0$, $\frac{z^2}{|z^2|} \frac{1}{\sqrt{1 + |z|}} \to e^{2i\theta }$. Thus, $(0, e^{2i\theta} ) \in \tg_0 X$ for all $\theta$. This implies that $\tg_0 X = \{z = 0\}$. 

Write $z = z_1 + i z_2$ and  $h = h_1 + i h_2$. For $x = (z,h(z))$, $z \neq 0$, $\tg_x X$ is generated by two vectors $u = (1, 0, \frac{\partial h_1}{\partial z_1}, \frac{\partial h_2}{\partial z_1})$ and $ v = (0, 1,  \frac{\partial h_2}{\partial z_1},\frac{\partial h_2}{\partial z_2})$. Denote by $\partial_1$ and $\partial_2$ the directional derivatives in the variable $z$ along $z_1$-axis and $z_2$-axis respectively. We know that
\begin{align*}
\partial_1 h & = \frac{\partial h_1}{\partial z_1}  + i \frac{\partial h_2}{\partial z_1} \\
\partial_2 h &= -i\frac{\partial h_1}{\partial z_2}  + \frac{\partial h_2}{\partial z_2}
\end{align*} 
Note that $h = \frac{z^2}{(z \bar{z})^{\frac{3}{4}}}$. Computation gives,

$$\partial_1 h = \frac{5}{4} \frac{z}{|z|^{\frac{3}{2}}} \partial_1 z - \frac{3}{4} \frac{z^3}{|z|^{\frac{7}{2}}} \partial_1 \bar{z}  = \frac{5}{4} \frac{z}{|z|^{\frac{3}{2}}} - \frac{3}{4} \frac{z^3}{|z|^{\frac{7}{2}}} $$

$$\partial_2 h = \frac{5}{4} \frac{z}{|z|^{\frac{3}{2}}} \partial_2 z - \frac{3}{4} \frac{z^3}{|z|^{\frac{7}{2}}} \partial_2 \bar{z}  = \frac{5}{4} \frac{z}{|z|^{\frac{3}{2}}} + \frac{3}{4} \frac{z^3}{|z|^{\frac{7}{2}}}. $$
If $z$ tends to $0$, then $|\partial_1 h|$ and $|\partial_2 h|$ tend to $\infty$. Therefore, 
\begin{align*}
\|u\| &= (1 + |\partial_1 h|^2)^\frac{1}{2} \to \infty, \text{ and}\\
\|v\| &= (1 + |\partial_2 h|^2)^\frac{1}{2} \to \infty. 
\end{align*}
Hence $\lim_{x \to 0} \tg_x X =\tg_0 X$. This implies (2).

Now we show (3). Denote by $\pi_0: \bb R^4 \to \tg_0 X$ the orthogonal projection from $\bb R^4$ onto $\tg_0 X$ : 
$$\pi_0|_X: X \ni (z, h(z)) \mapsto  h(z) \in \tg_0 X.$$
This map is not injective in any neighborhood of $0$ since $\forall z,\; h(z) = h(-z)$. By Theorem  \ref{thm_gluck}, $X$ is not a $C^1$ manifold. Since $X$ is a connected topological manifold, it cannot be the union of two or more $C^1$ manifolds of dimension $2$.  
\end{proof}
\end{exa}

\begin{dfn}[\cite{comte},\cite{co-me},\cite{ku-ra}]\label{def_density}\rm
Let $X\subset \bb R^n$ be a definable set and let $x \in \bb R^n$. Suppose that  $\dim X = k$. It is known that the following limit always exists
$$\theta(X, x) = \lim_{r\to 0}  \frac{\al H^k(X \cap \overline{\bff B}^n(x, r))}{r^k \mu_k}.$$
We call it the\textit{ density }of $X$ at the point $x$. 
\end{dfn}

\begin{rem} The notions of lower density and density are not the same even for definable sets. For example, consider $X: = \{(x, y, z) \in \bb R^3: z^4 = x^2 + y^2\}$. It is easy to see that $\Theta(X, 0) = \infty$ while $\theta (X, 0) = 0$. However, if  $X$ is a definable set and $\tg X$ is a trivial bundle then $\dim_{\al H} \tg_x X = \dim \tg_x X = \dim X$, and therefore $\Theta(X, x) = \theta(X, x)$ for every $x \in X$.
\end{rem}

\begin{lem}\label{rem_density} 
\begin{enumerate} 
\item Let $X, Y$ be definable sets of the same dimension. If $\dim (X \cap  Y) < \dim X$, then $\theta( X\cup Y, x) = \theta(X,x) + \theta (Y, x)$. If $X\subset Y$, then  $\theta( X, x) \leq \theta(Y,x)$.

\item If $X$ is a definable set then $\theta(X, x) \geq \theta(\tg_x X, 0)$.
\end{enumerate}
\end{lem}

\begin{proof}
$(i)$ is a direct consequence of the definition of density, and $(ii)$ of Theorem 3.8 of \cite{ku-ra}.
\end{proof}

\begin{thm}\label{thm_gh_ho_general} Let $X$ be a locally closed definable subset of $\bb R^n$. If  $\tg X$ is a continuous trivial bundle and $\theta(X, x) < 3/2$ for every $x \in X$, then $X$ is a $C^1$ manifold. 
\end{thm}

\begin{rem}The condition on the density above is sharp, since $$X:=\{(x,y);\; y=0\}\cup\{(x,y);\; x>0, y=x^2\}$$ satisfies all other hypothesis of Theorem \ref{thm_gh_ho_general} and $\theta(X,0)=3/2$.\end{rem}
\begin{proof}
Denote by $\pi_x: \bb R^n \to  \tg_x X$ the orthogonal projection. Suppose on the contrary that $X$ is not a $C^1$ manifold.  By Theorem \ref{thm_gluck}, there exists $x \in X$ such that there is no neighborhood of $x$ in $X$ to which the restriction of $\pi_x$ is injective. We may assume that $x$ coincides with the origin $0$ and $\tg_x X= \bb R^k$ where $k = \dim X$. The map $\pi_x$ now becomes $\pi: \bb R^n \to \bb R^k$, the orthogonal projection to the first $k$ coordinates. 

By Lemma \ref{lem_open_map}, there is an open neighborhood $U$ of $0$ in $X$ such that $\pi|_ U$ is an open map. Hence there exists $r > 0$ such that $\bff B^{k}(0,r) \subset \pi (U)$. Shrinking $U$ if necessary we assume that $\pi(U) = \bff B^{k}(0,r)$. We also assume that for every $x \in U$, $\tg_x X$ is not orthogonal to $\bb R^k$. 

By the uniform finiteness on fibres of definable sets, there exists $N \in \bb N$ such that for every $z \in\bff B^{k}(0,r)$, $|\pi^{-1} (z) \cap U|$, the number of connected components of $\pi^{-1} (z) \cap U$, does not exceed $N$. In fact, in this case, $|\pi^{-1} (z) \cap U| = \card (\pi^{-1} (z) \cap U)$ where card denotes the cardinality. If otherwise, there is a connected component of $\pi^{-1} (z) \cap U$, write $F$, such that $\dim F \geq 1$. Since $F$ is definable, there is a point $\tilde{z} \in K$ such that $\dim \tg_{\tilde{z}} F = \dim F \geq 1$. But $F \subset (\{z\}\times \bb R^{n-k})$, so $\tg_{\tilde{z}}F \subset (\{0\}\times \bb R^{n-k})$. This implies that $\tg_{\tilde{z}}F \subset \tg_{\tilde{z}} X  \cap \{0\}\times \bb R^{n-k}$. Since $\tilde{z} \in U$,  $\tg_{\tilde{z}} X$ is not orthogonal to $\bb R^k$, then $ \tg_{\tilde{z}} X  \cap (\{0\}\times \bb R^{n-k}) = \{0\}$, which gives a contradiction. 

Set
$$S_\kappa: = \{ z\in \bff B^{k}(0,r): \card(\pi^{-1} (z) \cap U )= \kappa\}.$$
Then $\{S_\kappa\}_{\kappa=1}^N$ becomes a definable partition of $ \bff B^{k}(0,r)$. We may assume that $0 \in \overline{S_\kappa}, \forall \kappa$ and  $S_N \neq \emptyset$. Note that $N \geq 2$ since the restriction of $\pi$ is not injective on any neighborhood of $0$. We claim that

(a) $S_N$ is an open set,

(b) $|(\pi^{-1}(S_N)\cap U| = N$,

(c) Each connected component of $\pi^{-1}(S_N)\cap U$ is a $C^1$ manifold. 

We now give a proof of the claim. 

Let $q \in S_N$. Since $\card(\pi^{-1}(q)\cap U) = N$, we may write $ \pi^{-1}(q)\cap U = \{q_1, \ldots, q_N \}$.  There is $\epsilon> 0$ sufficiently small such that $ K_i \cap K_j = \emptyset$ for $i \neq j$, where $K_i : = \bff B^n(q_i, \epsilon) \cap U$, $i \in \{1, \ldots, N\}$. Since the map $\pi|_U$ is open, there exists an open neighborhood $V_q$ of $q$ in $\bb R^k$ such that $V_q \subset \pi(K_i), \forall i \in \{1, \ldots, N\}$.  For $q' \in V_q$, $\forall i \in \{1, \ldots, N\}$, $\pi^{-1}(q') \cap K_i \neq \emptyset$, hence $\card(\pi^{-1}(q') \cap U) \geq N$. By the definition of $N$, $\card(\pi^{-1}(q') \cap U)$ cannot exceed $N$, therefore $\card(\pi^{-1}(q') \cap U) = N$, meaning $q'\in S_N$. This implies $V_q \subset S_N$, or equivalently $S_N$ is open, (a) is proved.

Denote by $A_1, \ldots, A_m$  the connected components of $\pi^{-1}(S_N) \cap U$. Note that for every $i$, $A_i$ is an open set in $U$, so $\pi(A_i)$ is an open set in $ \bb R^k$ since $\pi|_U$ is an open map. Assume that (b) does not hold, i.e., $m > N$. Then, there exists an $i \leq m$ such that $\pi(A_i)$ does not cover the whole of $S_N$, for simplicity we assume $i = 1$. Since $\pi(A_1) \subsetneq S_N$ is open, there exists  $p \in S_N \cap \partial \pi(A_1)$. Writing $\pi^{-1}(p) \cap U = \{p_1, \ldots, p_N\}$, there is an $\alpha \in \{1, \ldots, N\}$ such that $p_\alpha$ belongs to $\overline{A_1}$. But $p_\alpha \not \in A_1$ since $\pi(p_\alpha) \in \partial \pi(A_1)$, then  $p_\alpha \in A_\beta, \beta \neq 1$. This impiles that  $\overline{A_1} \cap A_\beta \neq \emptyset$, so $A_1$ and $A_\beta$ are the same connected component, which is a contradiction.

It follows from (b) that for each $i \in \{1, \ldots, N\}$ the restriction $\pi|_{A_i}: A_i \to S_N$ is a bijection. In other words, $A_i = \gr_{\xi_i}$ with $\xi_i: S_N \to \bb R^n$, $\xi_i(y) = \pi^{-1}(y) \cap A_i$. Note that $\tg A_i = \tg X|_{A_i}$ which is a continuous trivial bundle, its fibers are, moreover, not orthogonal to $\bb R^k$ by the construction. This shows that the function $\xi_i$ satisfies conditions of Lemma \ref{lem_C^1_map}, hence $\gr_{\xi_i}$ is a $C^1$ manifold. This ends the proof of (c).

Let $z \in S_N$ such that  $\|z\| < r/4$. Let $z'$ be a point realizing the distance from $z$ to the boundary $\partial S_N$ of $S_N$. Since $0 \in \partial S_N$, $s := \|z' - z\| \leq r/4$, and then
$\|z'\| \leq \|z' - z\| + \|z\| \leq r/4 + r/4 = r/2$, hence $z' \in \pi (U) = \bff B^k(0, r)$. Since $\bff B^k(z, s) \subset S_N$,  $\pi^{-1}(\bff B(z, s))\cap U$ has exactly $N$ connected components, denoted by $\{B_1, \ldots, B_N\}$. Remark that $z' \not \in S_N$ (i.e, $\card(\pi^{-1}(z') \cap U) < N$) but $z' \in \pi(\overline{B_i}), \forall i \in \{1, \ldots, N\}$, so there are $i, j$, $i\neq j$ such that $\overline{B_i} \cap \overline{B_j} \cap \pi^{-1}(z') \neq \emptyset$. Take $w \in \overline{B_i} \cap \overline{B_j} \cap \pi^{-1}(z')$. Take $C\subset \bff B^k(0, r)$ a small closed ball outside $\overline{\bff B}^k(z,s)$ and tangent to $\overline{\bff B}^k(z,s)$ at $z'$. Denote by $D$ the connected component of $\pi^{-1}(C) \cap U$ which contains $w$.

Since $\tg_w U$ is a $k$-dimensional linear subspace of $\bb R^n$ and $\pi|_{\tg_w U}: \tg_w U \to \bb R^k$ is  a linear bijective map, $D, B_i, B_j$ are disjoint definable sets of dimension $k$ in $U$ and $\tg_w D$, $\tg_w B_i$, $\tg_w B_j$ are half $k$-planes. By Lemma \ref{rem_density},
\begin{align*}
\theta (X, w) = \theta (U, w) &\geq \theta (B_i, w) + \theta (B_j, w) +\theta (D,w)\\
& \geq  \theta( \tg_wB_i, 0) +\theta(\tg_w B_j, 0) + \theta(\tg_w D, 0) \\
 & = \frac{1}{2} + \frac{1}{2} + \frac{1}{2 } = \frac{3}{2}.
\end{align*} 
This contradicts the hypothesis that $\theta (X, x) < \frac{3}{2}$ for every $x \in X$. 
\end{proof}

\subsection*{Acknowledgements} This research has been  supported by ANR project STAAVF. The third author has been also received the funding from NCN grant 2014/13/B/ST1/00543.

\end{document}